\documentclass[11pt]{article}

\usepackage{amsmath,amssymb,amsfonts,amsthm}

\newcommand{\Ann}{\mbox{Ann}\,}

\newcommand{\Ext}{\mbox{Ext}\,}
\newcommand{\Tor}{\mbox{Tor}\,}
\newcommand{\Spec}{\mbox{Spec}\,}
\newcommand{\codim}{\mbox{codim}\,}

\newcommand{\Ker}{\mbox{Ker}\,}

\newcommand{\Supp}{\mbox{Supp}\,}
\newcommand{\gr}{\mbox{grade}\,}
\newcommand{\depth}{\mbox{depth}\,}
\renewcommand{\dim}{\mbox{dim}\,}
\newcommand{\cmd}{\mbox{cmd}\,}

\newcommand{\rank}{\mbox{rank}\,}
\newcommand{\cx}{\mbox{cx}}
\newcommand{\pd}{\mbox{pd}\,}

\newcommand{\gsd}{\mbox{G$^*-$dim}\,}
\newcommand{\qpd}{\mbox{qpd}\,}
\newcommand{\cd}{\mbox{CI-dim}\,}

\renewcommand{\H}{\mbox{H}}
\newcommand{\V}{\mbox{V}}

\newcommand{\T}{\mathrm}

\newcommand{\lo}{\longrightarrow}

\newcommand{\fm}{\frak{m}}
\newcommand{\fp}{\frak{p}}
\newcommand{\fq}{\frak{q}}

\date{}
\begin{document}

\title{\bf Special homological dimensions and Intersection Theorem \footnotetext{The second author is supported by a grant
from IPM (NO. 82130212).} }

\author{Tirdad Sharif and Siamak Yassemi}
\maketitle

\begin{abstract}

\noindent Let $(R,\fm)$ be commutative Noetherian local ring. It
is shown that $R$ is Cohen--Macaulay ring if there exists a
Cohen--Macaulay finite (i.e. finitely generated) $R$--module with
finite upper Gorenstein dimension. In addition, we show that, in
the Intersection Theorem, projective dimension can be replaced by
quasi--projective dimension.

\end{abstract}

\baselineskip=18pt

\vspace{.3in}

\noindent{\bf 1. Introduction}

\vspace{.2in}

Let $M$ and $N$ be finite $R$--modules and $\pd_RM<\infty$. The
New Intersection Theorem of Peskine and Szpiro [{\bf PS}],
Hochster [{\bf H}], and P. Roberts [{\bf R1}], [{\bf R2}] yields
an inequality
$$(1) \,\,\,\,\,\,\dim_RN\le\dim_R(M\otimes_RN)+\pd_RM.$$
By applying the inequality (1) with $N=R$ one derives,
$$(2) \,\,\,\,\,\,\dim R\le\dim_RM+\pd_RM.$$
By the Auslander--Buchsbaum Formula, the inequality (2) is
equivalent to
$$(3)\,\,\,\,\,\, \cmd R\le\cmd_RM$$ where $\cmd_RM=\dim_RM-\depth_RM$ is the {\it Cohen--Macaulay
defect} of $M$ that is a non--negative integer which determines
the failure of $M$ to be Cohen--Macaulay; we set $\cmd R=\cmd_RR$.

It is shown that the projective dimension can be replaced by the
quasi--projective dimension (cf. [{\bf AGP}]) in (1).

 Also we show that the projective dimension can be replaced by upper Gorenstein dimension (cf. [{\bf
V}]) in (2) and (3). In addition, it is shown that the grade of a
module of finite upper Gorenstein dimension is actually equal to
its codimension.

\vspace{.1in}

\noindent (1.1) {\bf Setup and Notions.} Throughout, the rings
will denote a non--trivial, commutative, local, Noetherian ring
and the modules are finite (that means finitely generated). A
diagram of local homomorphisms $R\rightarrow R'\leftarrow Q$, with
$R\rightarrow R'$ a flat extension and $R'=Q/(x)$ where $x=x_1,
x_2, \cdots, x_c$ is a Q--regular sequence, is called a
quasi--deformation of codimension $c$. The {\it quasi--projective
dimension} of the $R$--module $M$ is defined by Avramov, Gasharov,
and Peeva [{\bf AGP}] as
$$\qpd_RM=\inf\{\pd_Q(M\otimes_RR')|\,\,R\rightarrow R'\leftarrow
Q\,\,\, \mbox{is a quasi--deformation}\}.$$

A local surjection $\pi\colon Q\to R$ is a Gorenstein deformation
if $\Ker(\pi)$ is a Gorenstein ideal.
A Gorenstein quasi--deformation of $R$ is a diagram of local
homomorphisms $R\rightarrow R'\leftarrow Q$, with $R\rightarrow
R'$ a flat extension and $R'\leftarrow Q$ a Gorenstein
deformation. The {\it upper Gorenstein dimension} of the
$R$--module $M$ is defined by Veliche [{\bf V}] as
$$\gsd_RM=\inf\{\pd_Q(M\otimes_RR')-\pd_QR'|\,\,R\rightarrow R'\leftarrow
Q\,\,\, \mbox{is a Gorenstein quasi--deformation}\}.$$

\noindent The $n$th {\it Betti number} of $M$ over $R$ is defined
by $\beta^R_n(M)=\rank_k(\Ext_R^n(M,k))$. The {\it complexity} of
$M$ is defined by
$$\cx_R M=\inf\{d\in \Bbb N_0|\,\,\beta^R_n(M)\le an^{d-1}\,\,\,
\mbox{for some positive real $a$ and}\,\,\, n\gg 0\}.$$ By [{\bf
AGP}; Thm. 5.11] for any $R$--module $M$ there is an equality
$\qpd_RM=\cd_RM+\cx_RM$ where $\cd_RM$ is the complete
intersection dimension of $M$.

An $R$-complex $X$ is a sequence of $R$-modules $X_\ell$ and
$R$-linear maps $\partial_\ell^X, \ell\in \mathbb{Z}$,
$$X=\cdots \lo X_{\ell+1} \stackrel{\partial^X_{\ell+1}}{\lo} X_\ell
\stackrel{\partial_\ell^X}{\lo} X_{\ell-1} \lo \cdots$$ such that
$\partial_\ell^X \partial_{\ell+1}^X=0$ for all $\ell\in
\mathbb{Z}$. The module $X_\ell$ is called the module in degree
$\ell$, and the map $\partial_\ell^X:X_\ell\lo X_{\ell-1}$ is
called the $\ell$-th differential. An $R$--module $M$ is thought
of as the complex $M=0\lo M\lo 0.$ The {\em supremum} and {\em
infimum} of $X$ are defined by
\[ \begin{array}{rl} \sup X &=
\sup\{{\ell} \in \Bbb Z |\H_{ \ell }(X) \neq 0 \}
\\ [.1in]\inf X &= \inf\{ {\ell} \in \Bbb Z | \H_{
\ell }(X) \neq 0 \}
\end{array} \]
A morphism $\alpha: X\lo Y$ is said to be a quasi-isomorphism if
the induced morphism $\T{H}(\alpha):\T{H}(X)\lo \T{H}(Y)$ is an
isomorphism. The {\it derived category} ${\cal D}(R)$ of the
category of $R$--complexes is the category of $R$--complexes
localized at the class of all quasi--isomorphisms. The full
subcategory ${\cal D}_{b}^f (R)$ consist of complexes $X$ with
$\H_{\ell}(X)$ a finite $R$--module for all $\ell$ and
$\H_{\ell}(X)=0$ for $|\ell|\gg 0$. The left derived functor of
the tensor product functor of $R$-complexes is denoted by
$-\otimes_R^{\T{\mathbf L}}-$. For a complex $X$, the dimension of
$X$ is defined by
$$\dim_RX=\sup\{\dim R/\fp-\inf X_{\fp}|\fp\in\Spec R\}.$$
When $M$ is an $R$--module, this notion agrees with the usual
definition of $\dim_RM$.

\noindent Let $R$ be a local and $X$ be a homologically finite
complex of $R$--modules. the $\cd_RX$ is defined by
Sather-Wagstaff in [{\bf S}] as follows
$$\cd_RX=\inf\{\pd_Q(X\otimes_R R')-\pd_QR'|\,\,R\rightarrow R'\leftarrow
Q\,\,\, \mbox{is a quasi--deformation}\}.$$

\vspace{.3in}

\noindent{\bf 2. Upper Gorenstein dimension.}

\vspace{.2in}

In this section it is shown that the ring $R$ is Cohen--Macaulay
if there exists a Cohen--Macaulay finite $R$--module $M$ of finite
upper Gorenstein dimension. (The converse is easy.) In addition it
is shown that the grade of a module of finite upper Gorenstein
dimension is actually equal to its codimension.

\vspace{.1in}

\noindent (2.1) {\bf Theorem.} Let $M$ be a finite $R$--module
with finite upper Gorenstein dimension. Then the following hold

\begin{verse}

(a) $\cmd R\le\cmd_RM$.

(b) $\dim R\le\dim_RM+\gsd_RM.$

\end{verse}

\vspace{.1in}

\noindent{\it Proof.} (a) Since $\gsd M<\infty$, there exists a
quasi--deformation $R\rightarrow R'=Q/J\leftarrow Q$ with
$\pd_QM'<\infty$, where $M'=M\otimes_RR'$. By the Intersection
Theorem, $\cmd Q\le\cmd_Q M'$. It is well--known that
$\cmd_QM'=\cmd_{R'}M'$. Now since $R\rightarrow R'$ is a flat
extension, by [{\bf AF}; 1.2] the following hold
$$\cmd R'=\cmd R+\cmd_R R'/\fm R';$$
$$\cmd_{R'}M'=\cmd_RM+\cmd_R R'/\fm R'.$$
On the other hand by [{\bf AFH}; Cor. 3.12] we have $\cmd Q=\cmd
R'$. Therefore
\[ \begin{array}{rl}
\cmd R+\cmd_R R'/\fm R'&\, =\cmd Q\\
&\, \le\cmd_QM'\\
&\, =\cmd_RM+\cmd_R R'/\fm R'
\end{array} \]

(b) This part is obtained by applying (a) and the equality $\depth
R=\depth_RM+\gsd_RM$, cf. [{\bf V}; Prop. 2.4].\hfill$\square$

\vspace{.2in}

Let $M$ be a finite $R$--module. The {\it grade} of $M$, $\gr_RM$,
defined by Rees to be the maximal length of $R$--regular sequence
in the annihilator of $M$. Also, the {\it codimension} $\codim_RM$
of the support of $M$ in the spectrum of $R$ is defined as the
height of the annihilator of $M$. In [{\bf AF}] Avramov and Foxby
have studied some properties of the codimension of a module with
finite projective dimension. Now as an application of Theorem
(1.2), we give a generalization of [{\bf AF}; Prop. 2.5] for a
module with finite upper Gorenstein dimension. Often it is
convenient to compute $\gr_RM$ and $codim_RM$ from the formulas
$$\gr_RM=\inf\{\depth R_\fp|\fp\in\Supp_RM\};$$
$$\codim_RM=\inf\{\dim R_\fp|\fp\in\Supp_RM\}.$$

\vspace{.1in}

\noindent (2.2) {\bf Proposition.} (See [{\bf AF}; Prop. 2.5]) If
$M$ is a non--zero finite $R$--module of finite upper Gorenstein
dimension, then $\gr_RM=\codim_RM$ and there exists a prime ideal
$\fp$ minimal in $\Supp_RM$ such that $R_\fp$ is Cohen--Macaulay
of dimension $\gr_RM$.

\vspace{.1in}

\noindent{\it Proof.} Choose $\fq\in\Supp_RM$ such that
$\gr_RM=\depth R_\fq$, and then choose $\fp$ contained in $\fq$
and minimal in $\Supp_RM$. By using [{\bf V}; prop. 2.4] and [{\bf
V}; prop. 2.10] we conclude from the choices of $\fq$ and $\fp$
that
$$\gr_RM=\depth
R_\fq\ge\gsd_{R_\fq}M_\fq\ge\gsd_{R_\fp}M_\fp=\depth
R_\fp\ge\gr_RM.$$ Therefore $\gr_RM=\depth R_\fp$. Since $M_\fp$
is an $R_\fp$--module of finite length and of finite upper
Gorenstein dimension, the ring $R_\fp$ is Cohen--Macaulay by
(1.2). Now by the following inequalities
$$\gr_RM=\dim R_\fp\ge\codim_RM\ge\gr_RM,$$
the assertion holds.\hfill$\square$

\vspace{.3in}

\noindent{\bf 3. Quasi--projective dimension.}

\vspace{.2in}

The main result in this section is the Theorem 3.3 That is a
generalization of the Intersection Theorem.

\vspace{.1in}

\noindent (3.1) {\bf Definition.} For $X\in{\cal D}^f_b(R)$ the
quasi-projective dimension of $X$ is defined as
$$\qpd_RX=\inf\{\pd_Q(X\otimes_RR')|\,\,R\rightarrow R'\leftarrow
Q\,\,\, \mbox{is a quasi--deformation}\}.$$

\vspace{.2in}

\noindent (3.2) {\bf Lemma.} If $\pi : R\to S$ is a surjection
local homomorphism and $Y, Z\in{\cal D}^f_b(S)$, then
$$\dim_R(Y\otimes_R^{\T{\mathbf L}}Z)=\dim_S(Y\otimes_S^{\T{\mathbf
L}}Z).$$

\vspace{.1in}

\noindent{\it Proof.} Let $I=\ker(\pi)$. Then we have
$R/\pi^{-1}(\frak q)\cong S/\frak q$ for every prime ideal $\frak
q$ of $S$.

Claim 1. $\Supp_R(Y\otimes_R^{\T{\mathbf L}}Z)\subseteq\V(I)$. Fix
$\frak p\in\Supp_R(Y\otimes_R^{\T{\mathbf
L}}Z)=\Supp_R(Y)\cap\Supp_R(Z)$. Then $\frak p\in\Supp_R(\H_i(Y))$
for some $i$. Since $\H_i(Y)$ is an $S$--module, one has
$I\subseteq\Ann_R(\H_i(Y))$ and therefore $\frak
p\in\Supp_R(\H_i(Y))=\V(\Ann_R(\H_i(Y)))\subseteq\V(I)$.

Claim 2. For a prime ideal $\frak q$ of $S$, one has $\frak
q\in\Supp_S(Y\otimes_S^{\T{\mathbf L}}Z)$ if and only if
$\pi^{-1}(\frak q)\in\Supp_R(Y\otimes_R^{\T{\mathbf L}}Z)$. Let
$\frak p=\pi^{-1}(\frak q)$. First note that the complexes
$Y_{\frak q}$ and $Y_{\frak p}$ are isomorphic over $R_{\frak p}$.
In particular, $\inf(Y_{\frak p})=\inf(Y_{\frak q})$ and $\frak
p\in\Supp_R(Y)$ if and only if $\frak q\in\Supp_S(Y)$. As we just
argued, $\frak p\in\Supp_R(Y\otimes_R^{\T{\mathbf
L}}Z)=\Supp_R(Y)\cap\Supp_R(Z)$if and only if $\frak
q\in\Supp_S(Y\otimes_S^{\T{\mathbf
L}}Z)=\Supp_S(Y)\cap\Supp_S(Z)$. This proves the claim.

Claim 3. For a prime ideal $\frak q$ of $S$, let $\frak
p=\pi^{-1}(\frak q)$; then
$$
\inf((Y\otimes_R^{\T{\mathbf L}}Z)_{\frak
p})=\inf((Y\otimes_S^{\T{\mathbf L}}Z)_{\frak q}).$$

One has
\[ \begin{array}{rl}
\inf((Y\otimes_R^{\T{\mathbf L}}Z)_{\frak p})\, &=\inf(Y_{\frak
p}\otimes_{R_{\frak p}}^{\T{\mathbf L}}Z_{\frak p})\\
\, &=\inf(Y_{\frak p})+\inf(Z_{\frak p})\\
\, &=\inf(Y_{\frak q})+\inf(Z_{\frak q})\\
\, &=\inf(Y_{\frak q}\otimes_{S_{\frak q}}^{\T{\mathbf L}}Z_{\frak
q})\\
\, &= \inf((Y\otimes_S^{\T{\mathbf L}}Z)_{\frak q}).
\end{array} \]

Now from claims 1--3 we have
\[ \begin{array}{rl}
\dim_R(Y\otimes_R^{\T{\mathbf L}}Z)\, &=\sup\{\dim(R/\frak
p)-\inf((Y\otimes_R^{\T{\mathbf L}}Z)_{\frak p})|\frak
p\in\Supp_R(Y\otimes_R^{\T{\mathbf L}}Z)\}\\
\, &=\sup\{\dim(S/\frak q)-\inf((Y\otimes_S^{\T{\mathbf
L}}Z)_{\frak q})|\frak
q\in\Supp_S(Y\otimes_S^{\T{\mathbf L}}Z)\}\\
\, &=\dim_S(Y\otimes_S^{\T{\mathbf L}}Z).
\end{array} \]
and this is the desired equality.\hfill$\square$

\vspace{.2in}

\noindent (3.3) {\bf Theorem.} Let $Y\in{\cal D}^f_b(R)$ with
finite quasi--projective dimension. Then for $X\in{\cal
D}^f_b(R)$,
$$\dim_R X\le \dim_R(Y\otimes_R^{\T{\mathbf L}}X)+\qpd_RY.$$

\vspace{.1in}

\noindent{\it Proof.} Since $\qpd_RY<\infty$, so there is a
quasi--deformation $R\rightarrow R'\leftarrow Q$ with
$\pd_QY'<\infty$ where $Y'=Y\otimes_RR'$. By the Intersection
Theorem (finite version) for complexes, cf. [{\bf F}; 18.5],
$$\dim_QX'\le\dim_Q(Y'\otimes_Q^{\T{\mathbf L}}X')+\pd_QY'.$$

We have the following
\[ \begin{array}{rl}
\dim_Q(Y'\otimes_Q^{\T{\mathbf L}}X')\, &
=\dim_{R'}(Y'\otimes_{R'}^{\T{\mathbf L}}X')\\
\, &=\dim_{R'}(Y\otimes_{R}^{\T{\mathbf L}}X)'\\
\, &=\dim_{R}(Y\otimes_{R}^{\T{\mathbf L}}X)+\dim R'/\frak m R'
\end{array} \]
\noindent where the the first equality comes from lemma 3.2. It is
also easy to see that
$$\dim_QX'=\dim_{R'}X'=\dim_RY+\dim R'/\frak m R'.$$
Thus the proof is completed. \hfill$\square$

\vspace{.2in}

In the Theorem (3.3), if we put finite $R$--modules $M$ and $N$
instead of the complexes $Y$ and $X$ then we get the following
result.

\noindent (3.4) {\bf Corollary.} Let $M$ be a finite $R$--module
with finite complete intersection dimension. Then for a finite
$R$--module $N$; $$\dim N\le\dim(N\otimes_RM)+\qpd_RM.$$

\vspace{.1in}

\noindent{\it Proof.} By [{\bf F}; 16.22 ],
$\dim(N\otimes_R^{\T{\mathbf L}}M)=\dim(N\otimes_RM)$. Therefore
the assertion is obtained by applying Theorem 3.3.\hfill$\square$

\vspace{.2in}

The following example shows that we can not replace
quasi-projective dimension with complete intersection dimension in
Theorem 3.3.

\vspace{.1in}

\noindent{\it Example} Let
$Q=k[|x_1,x_2,\cdots,x_n,y_1,y_2,\cdots,y_n|]$ where $k$ is a
field. Let $z_i=x_iy_i$ for $i=1,2,\cdots,n$. Consider the
$Q$--ideals $I=(z_1,z_2,\cdots,z_n)$, $J=(x_1,x_2,\cdots,x_n)$,
and $L=(y_1,y_2,\cdots,y_n)$. Let $R=Q/I$. Then $R$ is complete
intersection but is not regular. Consider the $R$--modules $A$ and
$B$ as $A=R/JR$ and $B=R/LR$. Then $A\otimes_RB=R/(J+L)R$ and
hence $\dim A\otimes_RB=0$. Since $R$ is complete intersection we
have $\cd_R A<\infty$ and hence, $\cd_R A=\depth R-\depth_RA=0$.
On the other hand $\dim B=n$ so $\dim B>\cd_RA+\dim A\otimes_RB$.

\vspace{.2in}

Let $M$ and $N$ be $R$--modules. Then we define
$$ \gr_R(M,N)=\inf\{i|\Ext^i_R(M,N)\ne 0\}. $$ \noindent If
$\Ext^i_R(M,N)=0$ for all $i$, then $\gr_R(M,N)=\infty$. In [{\bf
AY}; Theorem 3.1] Araya and Yoshino have used the Intersection
Theorem to prove ``Let $M$ and $N$ be finite $R$--modules with
$\pd_R N<\infty$ and $\Tor_i^R(M,N)=0$ for all $i>0$ then for any
finite $R$--module $L$, we have the following inequalities
$$
(3)\,\,\,\,\,  \gr_R(L,M)-\pd_R
N\le\gr_R(L,M\otimes_RN)\le\gr_R(L,M).''
$$
In the following Theorem we show that, in (3), projective
dimension can be replaced by quasi--projective dimension.

\vspace{.1in} \noindent (3.5) {\bf Theorem.} Let $M$ and $N$ be
finite $R$--modules with $\qpd_RN<\infty$ and $\T{Tor}^R_i(M,N)=0$
for any $i>0$. Then

$$\gr_R(L,M)-\qpd_RN\le\gr_R(L,M\otimes_RN)-\cx_RN\le\gr_R(L,M).
$$

\vspace{.1in}

\noindent{\it proof.} By [{\bf BH}; Prop. 1.2.10] there exists
$\fp\in\Supp_RL\cap\Supp_R(M\otimes_RN)$ such that
$$\gr_R(L,M\otimes_RN)=\depth_{R_\fp}(M_{\fp}\otimes_{R_\fp}N_{\fp}).$$ By [{\bf AGP}; Prop. 1.6]
$\cd_{R_\fp}N_\fp<\infty$ and so by applying [{\bf I}; Thm. 4.3]
we have
$$\gr_R(L,M\otimes_RN)=\depth_{R_\fp}N_\fp+\depth_{R_\fp}M_\fp-\depth_{R_\fp}R_\fp.$$ On the other hand
$\depth R_\fp-\depth_{R_\fp}N_{\fp}=\cd_{R_\fp}N_\fp$, see [{\bf
AGP}; Thm. 1.4]. Thus
$\gr_R(L,M\otimes_RN)=\depth_{R_\fp}M_\fp-\cd_{R_\fp}N_\fp$.
Therefore by applying [{\bf BH}; Prop. 1.2.10]
$$\gr_R(L,M\otimes_RN)\ge\gr_R(L,M)-\cd_RN.$$
Now the left inequality is obtained by applying the equality
$\cd_RN=\qpd_RN-\cx_RN$.

For the right inequality, there exists $\fp\in\Supp M\cap\Supp L$
such that $\gr_R(L,M)=\depth_{R_\fp}M_\fp$. Let $\fq$ be a minimal
element of the set $\Supp(R/\fp\otimes_RM\otimes_R N)$. Then
$\fp\subseteq\fq$ and so $\fq\in\Supp L$. By Corollary 3.2,
$\dim_{R_\fq}(R_\fq/\fp R_\fq\otimes_{R_\fq}M_\fq)\le
\qpd_{R_\fq}N_\fq$ and hence by [{\bf AGP}; Thm 5.11]
$\dim_{R_\fq}(R_\fq/\fp
R_\fq\otimes_{R_\fq}M_\fq)\le\cd_{R_\fq}N_{\fq}+\cx_{R_\fq}N_{\fq}$.
The following inequalities hold:
\[ \begin{array}{rl}
\cd_{R_\fq} N_\fq &\, \ge\dim_{R_\fq}(R_\fq/\fp R_\fq\otimes_{R_\fq}M_\fq)-\cx_{R_\fq}N_\fq\\
&\,\ge\dim_{R_\fp}(R_\fp/\fp R_\fp\otimes_{R_\fp}M_\fp)+\dim_{R_\fq}R_\fq/\fp R_\fq-\cx_{R_\fq}N_\fq\\
&\, \ge \depth_{R_\fq}M_\fq-\depth_{R_\fp}M_\fp-\cx_{R_\fq}N_\fq\\
&\, \ge \depth_{R_\fq}M_\fq-\depth_{R_\fp}M_\fp-\cx_R N.
\end{array} \]
Therefore
\[ \begin{array}{rl}
\cx_R N+\gr_R(L,M)&\, \ge\cx_R N+\depth_{R_\fp}M_\fp\\
&\, \ge\depth_{R_\fq}M_\fq-\cd_{R_\fq}N_\fq\\
&\, =\depth_{R_\fq}(M_\fq\otimes_{R_\fq}N_\fq)\\
&\, \ge\gr_R(L,M\otimes_R N)
\end{array} \]
Now the assertion holds.\hfill$\square$

\vspace{.2in}

The following example shows that the term ``$\cx_RN$'' is
necessary in the Theorem (3.5).

\vspace{.1in}

\noindent{\it Example.} Let $R=k[|X,Y|]/(XY)$, $N=R/yR$ where $x$
(resp. $y$) is image of $X$ (resp. $Y$) in $R$. Since $R$ is
complete intersection so $\cd_RN<\infty$. Since $\depth_RN=1$ we
have $\cd_RN=0$, the quasi--deformation can be chosen as $R=R'$
and $Q=k[|X,Y|]$. Then it is easy to see that $\cx_R(N)=1$. Set
$M=R$ and $L=R/xR$. Then $\gr_R(L,M)=0$ and
$\gr_R(L,M\otimes_RN)=1$.

\vspace{.3in}

\noindent {\large\bf Acknowledgments.} The authors would like to
thank Sean Sather-Wagstaff, University of Illinois, for his
invaluable comments, specially for the lemma 3.2.

\vspace{.3in}

\baselineskip=16pt

\begin{center}
\large {\bf References}
\end{center}
\vspace{.2in}

\begin{verse}


[Ap] D. Apassov, {\em Almost finite modules}, Comm. Algebra {\bf
27} (1999), 919--931.

[Av] L. L. Avramov, {\em Homological asymtotics of modules over
local rings}, in ``Commutative Algebra,'' Vol. 15, pp. 33--62,
MSRI, Berkeley, 1982; Springer-Verlag, New York, 1989.

[AGP] L. L. Avramov; V. N. Gasharov; I. V. Peeva, {\em Complete
intersection dimension} Inst. Hautes \'{E}tudes Sci. Publ. Math.
{\bf 86} (1997), 67--114.

[AF] L. L. Avramov; H.-B. Foxby, {\em Cohen-Macaulay properties of
ring homomorphisms} Adv. Math. {\bf 133} (1998), 54--95.

[AFB] L. L. Avramov; H.-B. Foxby; B. Herzog, {\em Structure of
local homomorphisms}, J. Algebra {\bf 164} (1994), 124--145.

[F] H.-B. Foxby, {\em  Hyperhomological algebra and commutative
algebra}, Notes in preparation.

[H] M. Hochster, {\em Topics in the homological theory of modules
over commutative rings}, CBMS Regional Conf. Ser. in Math., vol.
24, Amer. Math. Soc. Providence, RI, 1975.

[I] S. Iyengar, {\em Depth for complexes, and intersection
theorems}, Math. Z. {\bf 230} (1999), 545--567.


[PS] C. Peskine; L. Szpiro, {\em Syzygies et multiplicit\'{e}s},
C. R. Acad. Sci. Paris S\'{e}r. A {\bf 278} (1974), 1421--1424.

[R1] P. Roberts, {\em Two applications of dualizing complexes over
local rings}, Ann. Sci. \'{E}cole Norm. Sup. (4) {\bf 9} (1976),
103--106.

[R2] P. Roberts, {\em Le th\'{e}or\`{e}me d'intersection}, C. R.
Acad. Sci. Paris S\'{e}r. I Math. 304 (1987), 177--180

[V] O. Veliche, {\em Construction of modules with finite
homological dimensions}, J. Algebra {\bf 250} (2002), 427--449.

[Y] S. Yassemi, {\em G--dimension}, Math. Scand. {\bf 77} (1995),
161--174.

\end{verse}

{\bf \small {The Authors' Addresses}} \\[0.3 cm]
{\footnotesize Tirdad Sharif, Department of Mathematics, University of Tehran, Tehran, Iran.\\
E-mail address: {\tt sharif@ipm.ir} \\[0.2cm]
Siamak Yassemi, School of Mathematics, Institute for Studies in
Theoretical Physics and Mathematics, P.O. Box 19395-5746, Tehran,
Iran, and Department of Mathematics, University of Tehran, Tehran, Iran.\\
E-mail address: {\tt yassemi@ipm.ir}}

\end{document}